\documentclass[12pt,twoside]{article}
\usepackage{graphicx,graphics, amsfonts, amsbsy,amssymb, amsmath, cite, array,arydshln}

\allowdisplaybreaks

\def\bpsp{\begin{pspicture}}
\def\epsp{\end{pspicture}}

\newtheorem{theorem}{Theorem}[section]

\newtheorem{lemma}[theorem]{Lemma}
\newtheorem{corollary}[theorem]{Corollary}

\newtheorem{conjecture}{Conjecture}

\usepackage{tikz}
\usetikzlibrary{chains,positioning}
\usepackage{amsmath,tikz,pgfplots}
\usetikzlibrary{decorations.pathmorphing}

\def\qed{\nolinebreak\hfill\rule{.2cm}{.2cm}\par\addvspace{.5cm}}

\begin{document}
\title{On a conjecture of Laplacian energy of trees}
\author{ Hilal A. Ganie$^{a}$, Bilal A. Rather$ ^{b} $, S. Pirzada$ ^{c} $  \\
$^{a}$ {\em Department of School Education, JK Govt. Kashmir, India}\\
$^{b,c}${\em  Department of Mathematics, University of Kashmir, Srinagar, India}\\
$^a$hilahmad1119kt@gmail.com,~~ $^b$bilalahmadrr@gmail.com\\
$^c$pirzadasd@kashmiruniversity.ac.in}
\date{}

\pagestyle{myheadings} \markboth{Hilal, Bilal and Pirzada}{On a conjecture of Laplacian energy of trees }
\maketitle
\vskip 5mm
\noindent{\footnotesize \bf Abstract.}
Let $G$ be a simple graph with $n$ vertices, $m$ edges having Laplacian eigenvalues $\mu_1, \mu_2, \dots, \mu_{n-1},\mu_n=0$.  The Laplacian energy $LE(G)$ is defined as $LE(G)=\sum_{i=1}^{n}|\mu_i-\overline{d}|$, where $\overline{d}=\frac{2m}{n}$ is the average degree of $G$.  Radenkovi\'{c} and Gutman conjectured that among all trees of order $n$, the path graph $P_n$ has the smallest Laplacian energy. Let $ \mathcal{T}_{n}(d) $ be the family of trees of order $n$ having diameter $ d $. In this paper, we show that Laplacian energy of any tree $T\in \mathcal{T}_{n}(4)$ is greater than the Laplacian energy of $P_n$, thereby proving the conjecture for all trees of diameter $4$. We also show the truth of conjecture for all trees with number of non-pendent vertices at most $\frac{9n}{25}-2$. Further, we give some sufficient conditions for the conjecture to hold for a tree of order $n$.

\vskip 3mm



\section{Introduction}
Let $ G =(V,E)$ be a simple graph with vertex set $V(G)= \{v_{1},v_{2},\dots,v_{n}\} $ and edge set $ E(G)=\{e_{1},e_{2},\dots,e_{m}\} $. The parameter $n=|V(G)|$ is the order while  $m=|E(G)|$ is the size of the graph $G$.
 For a vertex $ v_{i}\in V(G)$ let $ d(v_{i})$ (or simply by $ d_{i}$, if underling graph $ G $ is clear from the context) denotes the degree of vertex $ v_i $ in $ G .$ A vertex of degree one is called a pendent vertex or a leaf. A tree of order $n$ denoted by $ T_{n} $ is a connected acyclic graph (that is a graph having no cycle). A path $ P_{n} $ of order $ n $ is a tree with two pendent vertices. A star $ S_{n} $ on $ n $ vertices is a tree having $ n-1 $ pendent vertices and one vertex of degree $ n-1 $. For other notations and undefined terminology from graph theory and matrix theory, we refer reader to any standard text, such as \cite{DMCM,hj,sp}.

Let $ A(G) $  be the $ (0,1)$-adjacency matrix of $ G $ with eigenvalues $ \lambda_{1}(G)\geq\lambda_{2}(G)\geq\dots\geq\lambda_{n}(G) $. The concept of energy $\mathbb{E}(G)$ of a graph $G$ was introduced in 1978 by I. Gutman \cite{gut} and is defined as
\begin{align*}
\mathbb{E}(G)=\sum\limits_{i=1}^{n}|\lambda_{i}|.
\end{align*}
This quantity has been studied both chemically and mathematically. The concept of energy was not only extended to other type of matrices associated with a graph, in fact it has been extended to general matrices \cite{niki}. A lot of work has been done in this direction and the research is still on as there are various unsettled problems. For some recent work regarding graph energy and related result we refer to \cite{mh,phbr} and the references therein.

Let $ D(G)=diag(d_{1},d_{2},\dots,d_{n}) $ be the diagonal matrix of vertex degrees $ d_{i}=d_{G}(v_{i}), i=1,2,\dots,n$ of graph $ G. $ The  real symmetric positive semi-definite matrix  $ L(G)=D(G)-A(G) $ is the Laplacian matrix of the graph $G$ and its eigenvalues which we denote by  $ \mu_{1}(G)\geq\mu_{2}(G)\geq\dots\geq \mu_{n-1}\ge \mu_{n}(G)=0 $ are the Laplacian eigenvalues of $G$. The Laplacian energy $LE(G)$ of a graph $G$ was introduced in 2006 by Gutman and Zhou \cite{gutzhou} and  is defined as
\begin{align*}
LE(G)=\sum\limits_{i=1}^{n}\left |\mu_{i}-\overline{d}\right |,
\end{align*}
where $ \overline{d}=\frac{2m}{n} $ is the average degree of $ G $. Using the fact $\sum\limits_{1=i}^{n-1}\mu_{i}=2m$, it is easy to see that \cite{fritscher}
\begin{equation}\label{laplacian energy}
 LE(G)=2\left (\sum\limits_{i=1}^{\sigma}\mu_{i}-\sigma \overline{d}\right )=2\max_{1\leq k\leq n}\left( \sum\limits_{i=1}^{k}\mu_{i}-k\overline{d} \right ),
\end{equation}
where $\sigma$ is the number of Laplacian eigenvalues greater than or equal to average degree $ \overline{d}$. This quantity is extensively studied because of its applications and the connections with other important spectral parameters of the graph $G$. The parameter $\sigma$ is also an active component of the present research and some works mostly on trees can be found in the literature \cite{sin}. For some recent works on Laplacian energy and related results we refer to \cite{hpbt} and the references therein.

An important problem in the spectral graph theory is to determine the extremal value of some spectral parameter and the characterization of the graph/graphs which attain this  extremal  value among all the graph of order $n$ or among a class of graph of order $n$. This problem has been considered for energy $\mathbb{E}(G)$ of a graph  for various families (like trees, unicyclic graphs, etc) and the graphs attaining the maximum/minimum values are completely determined. For trees, it is  shown in \cite{gut2} that the star graph $S_n$ has the minimum energy while  the path graph $P_n$ has the maximum energy among all trees of order $n$.

Like energy of a graph, a natural extremal question involving the Laplacian energy of a graph is to ask for the extremal value of this parameter among all $n$-vertex graphs, and for the $n$-vertex graph (or for the family of such graphs)
whose Laplacian energy achieves this extremal value. However, this seems to be a hard problem, which,
to the best of our knowledge, is still open.  Radenkovi\'{c} and Gutman \cite{radenkovic} studied the correlation between the energy and the Laplacian energy of trees. They computed the energy and Laplacian energy of all trees up to $14$ vertices and formulated the following conjecture.

\begin{conjecture}:\label{conjecture}
Let $ T $ be a  tree on $ n $ vertices. Then
\begin{align*}
LE(P_{n})\leq LE(T) \leq LE(S_{n}).
\end{align*}
 \end{conjecture}

 For trees of diameter $3$, it is shown in \cite{trevisan} that the Conjecture \ref{conjecture} is true. In the same paper \cite{trevisan} the authors have shown  by direct computations that the Conjecture \ref{conjecture} is true for all trees up to $18$ vertices. The right inequality of the Conjecture \ref{conjecture} is shown to be true by Fritscher et al. in \cite{fritscher} for all trees of order $n$. For the left inequality, Chang et al. in \cite{chang} verified the Conjecture \ref{conjecture} for trees of diameter $4$ and $5$ with perfect matching. Recently, Rahman et al. \cite{rahman} considered some families of trees of diameter $4$ (not considered in \cite{chang}) and verified the truth of left inequality of the Conjecture \ref{conjecture} for trees belonging to these families. But as a whole the left hand inequality of the Conjecture \ref{conjecture} is still open. It is clear from (\ref{laplacian energy}) that if we want to compare the Laplacian energies of two graphs with same number of vertices and edges, we need information (either exact or very close to exact) about the sum of their largest Laplacian eigenvalues and about the number of Laplacian eigenvalues that are larger than the average degree. Both these informations are closely related to the structure of the graph and  so considering this problem among a class of graphs is a very hard problem.\\
 \indent In this paper our aim is to verify the truth of left hand inequality of the Conjecture \ref{conjecture} for trees of diameter $4$ completely.

It is worth to mention that Conjecture \ref{conjecture} is also of interest from Matrix theory point of view. Let $\mathbb{M}_n(\mathbb{C})$ be the set of all square matrices of order $n$ with complex entries. The trace norm of a matrix $M\in \mathbb{M}_n(\mathbb{C})$ is defined as $\lVert M\rVert_{*}=\sum\limits_{i=1}^{n}\sigma_i(M)$, where $\sigma_1(M)\ge \sigma_2(M)\ge\cdots\ge \sigma_n(M)$ are the singular values of $M$(i.e. the square roots of the eigenvalues of $MM^*$, where $M^*$ is the complex conjugate of $M$). It is well known that for a symmetric matrix $M$, if $\sigma_i(M)$ are the singular values and  $\lambda_i(M)$, $ i=1,2,\dots,n,$  are the eigenvalues, then $ \sigma_i(M)=|\lambda_i(M)|$. In the light of this definition, it follows that the Laplacian energy $LE(G)$ of a graph $G$ is the trace norm of the matrix $ L(G)-\frac{2m}{n}I_n$, where $I_n$ is the identity matrix of order $n$. It is an interesting problem in Matrix theory to determine among a given  class of matrices the matrix (or the matrices) which attain the maximum value and the minimum value for the trace norm. The trace norm of matrices associated with the graphs and digraphs are extensively studied. For some recent papers in this direction see  \cite{mr} and the references therein.
 Using trace norm concept we can reformulate  the Conjecture \ref{conjecture} as:
 \begin{align*}
 \lVert L(P_n)-\frac{2(n-1)}{n}I_n \rVert_{*}\le \lVert L(T)-\frac{2(n-1)}{n}I_n \rVert_{*}\le \lVert L(S_n)-\frac{2(n-1)}{n}I_n \rVert_{*}.
 \end{align*}
This gives another motivation for the study of left inequality of Conjecture \ref{conjecture}.

We use the computational software's, \textbf{Wolfram Mathematica} and \textbf{Sage Math}, for various computations in the paper. For more literature in this direction, we refer the following references \cite{sin,zhou} and the references therein.

The paper is organised as follows. In Section $ 2 $, we state some  already known results. In Section $ 3 $ we obtain a lower for the Laplacian energy of a tree $T$ and using this lower bound we obtain condition on the number of non-pendent vertices, which guarantees the  truth of Conjecture \ref{conjecture}. We also devise a recursive method which is helpful to verify the truth of the  Conjecture \ref{conjecture}. In Section $4$,  we verify truth of the conjecture \ref{conjecture} for the trees of diameter at most $ 4 $. Lastly, in Section $ 5 $, we give some sufficient conditions for the Conjecture \ref{conjecture} to hold for a tree of order $n$.

\section{Preliminaries}
In this section, we mention an algorithm and some known results which are helpful throughout the paper.

The following observation which gives the multiplicity of eigenvalue $1$ of the matrix $L(G)$ and can be found in \cite{faria}
\begin{lemma}\label{multiplicity of 1}\cite{faria}
Let $ p $ and $ q $ be the number of leaves and neighbors of $ G $, respectively. Then $ 1 $ is a Laplacian eigenvalue of $ G $ of multiplicity at least $ p-q\geq 0 $
\end{lemma}

The following lemma gives a lower bound for the $i^{th}$-largest Laplacian eigenvalue in terms of $i^{th}$-largest degree of the graph $G$.
\begin{lemma}\cite{brouwerheamers}\label{eigenvalue lower bound}
Let $\mu_{1}\geq \mu_{2}\geq\dots\geq \mu_{n} $ be the Laplacian eigenvalues of a connected graph $ G $. Then
\begin{align*}  \mu_{i}\geq d_{i}-i+2, \qquad i=1,2,\dots,n, \end{align*}
where $ d_{1}\geq\dots\geq d_{n} $ is the degree sequence of the vertices of $ G. $
\end{lemma}

The main tool to prove our results is the following algorithm  \cite{fricke,jacobs}, which can be used to
estimate the Laplacian eigenvalues of a given tree.\\
 \textbf{Algorithm (I).}\\
The algorithm associates with each vertex $ v $, a rational function $ a(v)=\frac{r}{s} $. Here $ r $ and $ s $ are members of the polynomial ring $ \mathbb{Q}[\lambda] .$ These are computed bottom-up starting with the leaves which are assigned $ \lambda-1 $ (the trees can be rooted in an arbitrary way). Once all the children of $ v $ have been processed, $ v $ is assigned the function
\begin{align*}
a(v)=\lambda-d_{v}-\sum_{c\in C} \frac{1}{a(c)},
\end{align*}
where $ C $ is the set of its children and $ d_{v} $ is the degree of $ v $. After all the vertices have been processed, we compute the characteristic polynomial by taking the product of all functions $ a(v) $:
\begin{align*}
p(\lambda)=\prod_{v\in V} a(v).
\end{align*}

\noindent\textbf{Algorithm (II)}:\cite{jacobs1} \\
\texttt{ Input: tree $T$, scalar $\alpha$}\\
\texttt{Output: diagonal matrix D congruent to $L(T)+\alpha I$}\\
Algorithm Diagonalize ($T,\alpha$)

initialize $a(v):=d(v)+\alpha,$ for all vertices $ v $

order vertices bottom up

\textbf{for} $ k=1\ \text{to}\ n $

\qquad \textbf{if} $ v_{k} $ is a leaf then continue

\qquad \textbf{else if} $ a(c)\neq 0 $ for all children $ c $ of $ v_{k} $ then

\qquad \qquad $ a(v_{k}):= a(v_{k})-\sum \frac{1}{a(c)},$ summing over all children of $ v_{k} $

\qquad\textbf{ else}

\qquad \qquad select one child $ v_{j} $ of $ v_{k} $ for which $ a(v_{j})=0 $

\qquad \qquad $ a(v_{k}):= -\frac{1}{2} $

\qquad\qquad $ a(v_{j}) := 2$

\qquad\qquad \textbf{if} $ v_{k} $ has a parent $ v_{l} $, remove the edge $ v_{k}v_{l} $

end loop

This algorithm is useful to calculate the number of eigenvalues of the Laplacian matrix of $T$ lying in a
given interval and hence helps us to estimate $\sigma$ for a tree. It is worth to mention that the diagonal elements of the output matrix correspond precisely to the values $a(v)$ on each node $v$ of the tree. The following observation due to Jacobs and Trevisan \cite{jacobs1} is helpful throughout the paper.
\begin{lemma}\cite{jacobs1}\label{number1}
 Let $ T $ be a tree and $ D $ be the diagonal matrix produced by the algorithm Diagonalize $ (T,-\alpha) $. Then the following assertions hold.\\
$ (a)  $ The number of positive entries in $ D $ is the number of Laplacian eigenvalues of $ T $ that are greater than $ \alpha. $\\
$ (b)  $ The number of negative entries in $ D $ is the number of  Laplacian eigenvalues of $ T $ that are smaller than $ \alpha. $\\
$ (c)  $ If there are $ j $ zero diagonal entries in $ D $, then $ \alpha $ is a Laplacian eigenvalues of $ T $ with multiplicity $ j $.
\end{lemma}

Using direct computation and the fact that
\begin{align*}
\frac{\pi}{n}\sum\limits_{j=1}^{\lfloor\frac{n}{2}\rfloor}\cos\frac{\pi j}{n}\le \int\limits_{0}^{\frac{\pi}{2}}\sin x dx=1,
\end{align*}
 the following upper bound for the Laplacian energy of a path $P_n$ was established in \cite{trevisan}.
\begin{lemma}\label{pathenergy}\cite{trevisan}
Let $ P_{n}  $ be a path on $ n $ vertices. Then
\begin{align*}
LE(P_{n})\leq 2+\frac{4n}{\pi}.
\end{align*}
\end{lemma}

The following observation can be found in \cite{groneme} and gives the interlacing property of the Laplacian eigenvalues of a graph and its spanning subgraph.
\begin{lemma}\label{lem1}
If $G^{\prime}=G+e$ is the graph obtained from $G$ by adding a new edge $e$, then the Laplacian
eigenvalues of $G$ interlace the Laplacian eigenvalues of $G^{\prime}$, that is,
\begin{align*}
\mu_1(G^{\prime})\geq\mu_1(G)\geq\mu_2(G^{\prime})\geq \mu_2(G)\geq\cdots\geq \mu_n(G^{\prime})\geq\mu_n(G)=0.
\end{align*}
\end{lemma}

 The following result about the distribution of Laplacian eigenvalues of trees can be found in \cite{sin}.
 \begin{lemma}\label{lem3}
 The number of Laplacian eigenvalues less than the average degree $2-\frac{2}{n}$ of a tree $T$ of order $n$ is at least $\lceil\frac{n}{2}\rceil$.
 \end{lemma}

\section{Laplacian energy of a tree}

In this section, we obtain a lower bound for the Laplacian energy of a tree $T$ and using this lower bound we obtain condition on the number of non-pendent vertices, which guarantees the  truth of Conjecture \ref{conjecture}. We also devise a recursive method which is helpful to verify the truth of Conjecture \ref{conjecture}.

For non-increasing real sequences $(x)=(x_1,x_2,\dots,x_n)$ and $(y)=(y_1,y_2,\dots,y_n)$ of length $n$, we say that $(x)$ is majorized by $(y)$ or $(y)$ majorizes $(x)$, denoted by $(x)\preceq (y)$ if
\begin{align*}
&\sum\limits_{i=1}^{n}x_i =\sum\limits_{i=1}^{n}y_i~~~\text{and}~~~\sum\limits_{i=1}^{k}x_i \leq\sum\limits_{i=1}^{k}y_i, ~~\text{for  all}~~k=1,2,\dots,n-1.
\end{align*}
 The eigenvalues of the Laplacian matrix $L(G)$ and majorization are connected by the
majorization relation between the eigenvalues and the diagonal elements. The following observation can be found in \cite{grone}.
\begin{lemma}\label{lem12}
Let $G$ be a  graph of order $n$ having Laplacian eigenvalues $\mu_1\ge \mu_2\ge \dots\ge \mu_n\ge 0$ and vertex degrees $d_{1}\ge d_{2}\ge \dots\ge d_{n}$. Then $(d_1,d_2,\dots,d_n)\preceq(\mu_1,\mu_2,\dots,\mu_n)$, that is
\begin{align*}
\sum\limits_{i=1}^{k}\mu_i \ge 1+\sum\limits_{i=1}^{k}d_i ~~\text{for  all}~~k=1,2,\dots,n-1.
\end{align*}
\end{lemma}

\begin{theorem}\label{thm1}
Let $T$ be a tree of order $n\ge 4$ and let $P_n$ be the path graph on $n$ vertices. If $T$ has $s$ internal vertices(non-pendent) vertices, then
\begin{align*}
LE(T)\ge LE(P_n),
\end{align*} provided that $\Big(\frac{\pi-2}{\pi}\Big)n\ge s+2-\frac{2s}{n}$.
\end{theorem}
\noindent{\bf Proof.} Let $\mu_1(T)\ge \mu_2(T)\ge \cdots\ge \mu_{n-1}(T)\ge \mu_n(T)=0$ be the Laplacian eigenvalues of $T$.  Let $s\ge 1$ be the number of internal vertices (that is, the vertices having degree greater than 1) of $T$, then the number of pendent vertices $p$ of $T$ is $p=n-s$. Since $1\le s\le n-2$, therefore using definition of Laplacian energy (\ref{laplacian energy}), we have
\begin{align}\label{u}
 LE(T)=&2\max_{1\leq k\leq n}\left( \sum\limits_{1=i}^{k}\mu_{i}(T)-k\overline{d} \right)\nonumber\\&
 \ge 2\Big( \sum\limits_{i=1}^{s}\mu_{i}(T)-s\overline{d} \Big).
\end{align}
Let $d_1(T)\ge d_2(T)\ge \cdots\ge d_n(T)$ be the vertex degree of $T$, then by Lemma \ref{lem12}, we have
\begin{align*}
\sum\limits_{i=1}^{k}\mu_{i}(T)\ge 1+\sum\limits_{i=1}^{k}d_i,~~~\text{for~~ all}~~~ k=1,2,\dots,n-1.
\end{align*}
Using this in (\ref{u}), we get
\begin{align}\label{u1}
LE(T)\ge 2\Big(1+\sum\limits_{1=i}^{s}d_i-s\overline{d} \Big).
\end{align}
Since $d_1+d_2+\cdots+d_n=2m=2n-2$ and $T$ has $p=n-s$ pendent vertices, it follows that $d_1+d_2+\cdots+d_s+n-s=2n-2$, which implies that $\sum\limits_{1=i}^{s}d_i=n+s-2$. Using this in (\ref{u1}), we get
\begin{align*}\label{u2}
LE(T)\ge 2n+2s-2-2s\overline{d}\ge 2+\frac{4n}{\pi},
\end{align*}
provided that $\Big(\frac{\pi-2}{\pi}\Big)n\ge s+2-\frac{2s}{n}$. Now, using Lemma \ref{pathenergy}, the result follow.\qed

The following observation is clear from the proof of Theorem \ref{thm1}.

\begin{corollary}\label{cor1}
For a tree $T$ on $n$ vertices having vertex degrees $d_1\ge d_2\ge\cdots\ge d_n$. If $k$, $1\le k\le n-1$ is any positive integer, then
\begin{align*}
LE(T)\ge 2\Big(1+\sum\limits_{1=i}^{k}d_i-k\overline{d} \Big).
\end{align*}
\end{corollary}
From this Corollary \ref{cor1}, it follows that any information about the degrees of a tree $T$ can be used to obtain a lower bound for Laplacian energy $LE(T)$ of $T$, which in turn can be helpful to get an insight for the possible solution  of the  Conjecture \ref{conjecture}.\\
\indent The following observation is also immediate from Theorem \ref{thm1}.
 \begin{corollary}\label{cor2}
 Let $T$ be a tree on $n$ vertices having Laplacian eigenvalues $\mu_1\ge \mu_2\ge\cdots\ge \mu_n$. If $k$, $1\le k\le n-1$ is any positive integer, then
 \begin{align*}
 LE(T)\ge 2\Big(\sum\limits_{i=1}^{k}\mu_i-k\overline{d} \Big).
 \end{align*}
 \end{corollary}
From  Corollary \ref{cor2}, it is clear that any lower bound for the sum of $k$ largest Laplacian eigenvalues  $S_k(T)=\sum\limits_{i=1}^{k}\mu_i$ of $T$  can be used to obtain a lower bound for Laplacian energy $LE(T)$, which in turn can be helpful to get an insight for the possible solution  of  Conjecture \ref{conjecture}.\\
The following observation gives the insight how Theorem \ref{thm1} can be helpful to prove the Conjecture \ref{conjecture}.
 \begin{corollary}\label{cor3}
Let $T$ be a tree of order $n\ge 4$ having $s$ internal vertices(non-pendent) vertices.\\
{\bf (i)} If $s=1$, then Conjecture \ref{conjecture} holds for all $n\ge 9$. {\bf (ii)} If $s=2$, then Conjecture \ref{conjecture} holds for all $n\ge 12$. {\bf (iii)} If $s=3$, then Conjecture \ref{conjecture} holds for all $n\ge 14$. {\bf (iv)} If $s=4$, then Conjecture \ref{conjecture} holds for all $n\ge 17$. {\bf (v)} If $s=5$, then Conjecture \ref{conjecture} holds for all $n\ge 20$. {\bf (vi)} If $s=6$, then Conjecture \ref{conjecture} holds for all $n\ge 23$. {\bf (vii)} If $s=7$, then Conjecture \ref{conjecture} holds for all $n\ge 25$. {\bf (viii)} If $s\le \frac{9n}{25}-2$, then Conjecture \ref{conjecture} holds for all $n$.
 \end{corollary}

The result gives a lower bound for the Laplacian energy of a tree $T$ in terms of sum of $k_i$ largest Laplacian eigenvalues of $T_i$, where $T_i$, for $i=1,2$ are components of $T$ obtained by deleting any non-pendent edge.
\begin{theorem}\label{thm2}
Let $T$ be a tree of order $n\ge 8$ and let $e$ be a non-pendent edge of $T$. Let $T-e=T_1\cup T_2$ and let $\sigma$ be the number of Laplacian eigenvalues of $T-e$ which are greater than or equal to average degree $\overline{d}(T-e)$. Then
\begin{align*}
LE(T)\ge 2S_{k_{1}}(T_1)+2S_{k_{2}}(T_2)-4\sigma+\frac{4\sigma}{n},
\end{align*} where $k_1,k_2$ are respectively, the number of Laplacian eigenvalues of $T_1$, $T_2$ which are greater than or equal to  $\overline{d}(T-e)$ with $k_1+k_2=\sigma$ and $S_k(T)$ is the sum of $k$ largest Laplacian eigenvalues of $T$.
\end{theorem}
\noindent{\bf Proof.} Let $\mu_1(T)\ge \mu_2(T)\ge \cdots\ge \mu_{n-1}(T)\ge \mu_n(T)=0$ be the Laplacian eigenvalues of $T$ having average degree $\overline{d}(T)=2-\frac{2}{n}$. Let $e$ be a non-pendent edge in $T$ and let $T_1$, $T_2$ be the components of $T-e$. Let $|V(T_i)|=n_i$, for $i=1,2$ and $\overline{d}(T-e)$ be the average degree of $T-e$. Then $n=n_1+n_2$ and $\overline{d}(T-e)=2-\frac{4}{n}$. Let $\sigma$ be the number of Laplacian eigenvalues of $T-e=T_1\cup T_2$, which are greater than or equal to average degree $\overline{d}(T-e)$, therefore $1\le \sigma\le n-2$.
It is well-known that $1\le \sigma\le n-2$, therefore using definition of Laplacian energy (\ref{laplacian energy}), we have
\begin{equation}\label{k}
 LE(T)=2\max_{1\leq k\leq n}\left( \sum\limits_{i=1}^{k}\mu_{i}(T)-k\overline{d}(T) \right)
 \ge 2\Big( \sum\limits_{i=1}^{\sigma}\mu_{i}(T)-\sigma(2-\frac{2}{n}) \Big).
\end{equation}
Let $\mu_1(T-e)\ge \mu_2(T-e)\ge \cdots\ge \mu_{n-2}(T-e)\ge\mu_{n-1}(T-e)= \mu_n(T-e)=0 $ be the Laplacian eigenvalues of $T-e$. By Lemma \ref{lem1}, we have $\mu_i(T)\ge \mu_i(T-e)$ giving that $\sum\limits_{i=1}^{\sigma}\mu_{i}(T) \ge \sum\limits_{i=1}^{\sigma}\mu_{i}(T-e)$. With this it follows from (\ref{k}) that
\begin{align}\label{v1}
 LE(T)\ge 2\sum\limits_{i=1}^{\sigma}\mu_{i}(T-e)-4\sigma+\frac{4\sigma}{n}.
\end{align}
Since Laplacian spectrum of $T_1\cup T_2$ is the union of Laplacian spectrum of $T_1$ and Laplacian spectrum of $T_2$. Let $k_i\ge 1$ be the number of Laplacian eigenvalues of $T_i$ which are greater than or equal to  $\overline{d}(T-e)$, then $k_1+k_2=\sigma$. Therefore, we have from (\ref{v1}) that
\begin{align*}
 LE(T)\ge 2\Big(\sum\limits_{i=1}^{k_1}\mu_{i}(T_1)+\sum\limits_{i=1}^{k_2}\mu_{i}(T_2)\Big)-4\sigma+\frac{4\sigma}{n}.
\end{align*}
The result now follows.\qed

Let $\sigma_i$ be the number of Laplacian eigenvalues of $T_i$ which are greater than or equal to average degree $\overline{d}(T_i)=2-\frac{2}{n_i}$. If $n_1\ge n_2$, the it is easy to see that $\overline{d}(T_1)\ge \overline{d}(T-e) $ and $\overline{d}(T_2)\le \overline{d}(T-e)$. Therefore, it follows that $k_1\ge \sigma_1$ and $k_2\le \sigma_2$. If $k_1=\sigma_1$ and $k_2=\sigma_2$, then we have the following observation.
\begin{corollary}\label{coru}
Let $T$ be a tree of order $n\ge 8$ and let $e$ be a non-pendent edge of $T$. Let $T_1$  and  $T_2$ be the components of $T-e$ and  let $\sigma$ be the number of Laplacian eigenvalues of $T-e$ which are greater than or equal to average degree $\overline{d}(T-e)$. Let  $k_i$ and $\sigma_i$ be  respectively, the number of Laplacian of $T_i$ which are greater than or equal to  $\overline{d}(T-e)$ with $k_1+k_2=\sigma$ and the number of Laplacian of $T_i$ which are greater than or equal to  $\overline{d}(T_i)$. If $\sigma_1=k_1$ and $\sigma_2=k_2$, then Conjecture \ref{conjecture} holds for $T$, provided that $LE(T_i)\ge 2+\frac{4n_i}{\pi}$ for $i=1,2$.
\end{corollary}
\noindent{\bf Proof.}
Let $\sigma_i$ be the number of Laplacian eigenvalues of $T_i$ which are greater than or equal to average degree $\overline{d}(T_i)=2-\frac{2}{n_i}$. Then, by definition of Laplacian energy (\ref{laplacian energy}), we have
\begin{align}\label{v2}
 LE(T_i)=2\sum\limits_{j=1}^{\sigma_i}\mu_{j}(T_i)-2\sigma_i \overline{d}(T_i).
\end{align}
Suppose that $LE(T_i)\ge 2+\frac{4n_i}{\pi}$ for $i=1,2$. Then, it follows from (\ref{v2}) that
\begin{align*}
2S_{\sigma_i}(T_i)\ge 2+\frac{4n_i}{\pi}+2\sigma_i (2-\frac{2}{n_i}).
\end{align*}
Now, if $\sigma_i=k_i$, for $i=1,2$, then it follows  from Theorem \ref{thm2} that
\begin{align*}
LE(T)\ge& 2+\frac{4n_1}{\pi}+2\sigma_1 (2-\frac{2}{n_1})+2+\frac{4n_2}{\pi}+2\sigma_2 (2-\frac{2}{n_2})-4\sigma+\frac{4\sigma}{n}\\&
=\frac{4n}{\pi}+4+\frac{4\sigma}{n}-\Big(\frac{4\sigma_1}{n_1}+\frac{4\sigma_2}{n_2}\Big)\\&
\ge \frac{4n}{\pi}+2,
\end{align*}
provided that $1+\frac{2\sigma}{n}\ge \frac{2\sigma_1}{n_1}+\frac{2\sigma_2}{n_2}$. Note that we have  used $\sigma_1+\sigma_2=k_1+k_2=\sigma$.  Since $n=n_1+n_2$ and $\sigma=\sigma_1+\sigma_2$, therefore $1+\frac{2\sigma}{n}\ge \frac{2\sigma_1}{n_1}+\frac{2\sigma_2}{n_2}$ gives that $n_1n_2(n_1+n_2)\ge 2\sigma_1n^{2}_2+2\sigma_2n^{2}_1$, which in turn gives that $n_1^{2}(n_2-2\sigma_2)+n_2^{2}(n_1-2\sigma_1)\ge 0$. Now applying Lemma \ref{lem3} to tree $T_i$, it follows that this last inequality always holds. This completes the proof.\qed

\section{Trees of diameter at most $ 4 $}

In this section, we verify the truth of Conjecture \ref{conjecture} for trees of diameter at most $4$.\\
 \indent Let $ \mathcal{T}_{n}(d)$ be the family of trees each  of diameter $d $ and order $ n\geq 3$. In particular, $ \mathcal{T}_{n}(4)$ is the family of trees with diameter $4$ and order $n$. In the family of trees $\mathcal{T}_{n}(4)$, there are more than $19$ possible subfamilies which we need to consider in order to verify the truth of Conjecture \ref{conjecture} for this family. Out of these subfamilies the truth of  Conjecture \ref{conjecture}  has already been verified for $5$ subfamilies; one in  \cite{chang} and four in \cite{rahman}.  Here we consider the family $\mathcal{T}_{n}(4)$ in general and verify the truth of Conjecture \ref{conjecture}.\\

Although Conjecture \ref{conjecture} has been already verified for tree of diameter at most $3$, see \cite{trevisan}. Here we prove it again by giving a very simple proof.
\begin{theorem}\label{thm23}
 Conjecture \ref{conjecture} holds for trees of diameter at most $3$.
\end{theorem}
\noindent{\bf Proof.} Let $T$ be  a tree of order $n\ge 3$ having diameter $2\le d\le 3$. If $d=2$, then $T$ is the star $K_{1,n-1}$ and so it follows by direct calculation that $LE(K_{1,n-1})\ge LE(P_n)$. If $d=3$, the $T$ is a double broom $T(a,b)$ of diameter $3$. Since the number of non-pendent vertices in $T(a,b)$ is $s=2$, it follows by part {\bf (ii)} of Corollary \ref{cor3} that  Conjecture \ref{conjecture} holds for all $n\ge 12$. For $3\le n\le 11$,  Conjecture \ref{conjecture} has been already verified by direct calculation. This completes the proof.\qed

Since Conjecture \ref{conjecture} is true for all trees with order up to $n\le 18$. Therefore, in the rest of this paper we consider $n\ge 19$. For the family $\mathcal{T}_n(4)$, the following result verifies Conjecture \ref{conjecture}.
\begin{theorem}\label{thm42}
 Conjecture \ref{conjecture} holds for the family $\mathcal{T}_n(4)$ of trees of diameter $4$.
\end{theorem}
\noindent{\bf Proof.} Let $T$ be a tree of order $n\ge 19$ having diameter $4$, then $T\in \mathcal{T}_n(4)$. Let $\mu_1(T)\ge \mu_2(T)\ge\cdots\ge\mu_{n-1}(T)\ge \mu_n=0$ be the Laplacian eigenvalues of $T$ and let $\overline{d}(T)=2-\frac{2}{n}$ be its average vertex degree. Using Lemma \ref{pathenergy}, it is clear that to show Conjecture \ref{conjecture} holds for $T$,  it suffices to show that the inequality
\begin{align}\label{a}
LE(T)\ge \frac{4n}{\pi}+2,
\end{align} holds for $T$, for all  $n\ge 19$.
Since diameter of $T$ is $4$, it follows that: \\
{\bf (i).} $T$ is either a double broom of diameter $4$ (shown in figure $1$) or;\\
{\bf (ii).} $T$ is a SNS-tree of diameter $4$ (shown in figure $3$).\\
\indent If $T$ is  a double broom of diameter $4$, then it clear from figure $1$ that $T$ has $s=3$ internal vertices. Therefore, using part {\bf (iii)} of Corollary \ref{cor3}, it follows that the inequality $LE(T)\ge \frac{4n}{\pi}+ 2$ always holds, giving that the result holds in this case. So, suppose that $T$ is a SNS-tree of order $n\ge 19$ having root vertex $v_0$ with $p\ge 0$ pendent vertices of level $0$, $r\ge 2$ vertices  of level $1$ such that each $v_i$ has $s_i$ pendent vertices attached, where at least two $s_i$ are non-zero. Clearly, order $n$ of $T$ in this case is $n=p+r+1+\sum\limits_{i=1}^{r}s_i$. If $p=0$ and  $s_i=1$, for all $i=1,2,\dots,r$, then $T$ is the tree $ T(4;2a,2b)$, $a+b=r$, shown in figure $1$. Applying algorithm $ {\bf (I)} $ to  $ T(4;2a,2b) $, we find that  the characteristic polynomial of $ T(4;2a,2b) $ and is given by $\phi(T(4;a,b),x)= x(x^2-3x+1)^{a+b-1}(x^2-x(a+b+3)+2a+2b+1).$
Applying algorithm $ {\bf (II)} $ with $\alpha=-1$ to the tree $ T(4;2a,2b) $ and using Lemma \ref{multiplicity of 1}, we see that $ a+b+1 $ eigenvalues are greater than one while as $ a+b $ eigenvalues are less than one. Since zeros of $ x^2-3x+1 $ and $x^2-x(a+b+3)+2a+2b+1  $ are $ \frac{3\pm\sqrt{5}}{2} $ and $ \frac{1}{2}\left (3+a+b\pm\sqrt{(a+b)^2-2(a+b)+5}\right ) $, respectively. It follows that $ \sigma= a+b$. Therefore, using $2(a+b)=n-1$, it follows  from the definition of Laplacian energy that
\begin{align*}
LE(T(4;a,b))=& 2\left (\sum\limits_{i=1}^{\sigma}\mu_{i}(T(4;a,b))-\sigma \overline{d}(T(4;a,b))\right )\\&
= 2\left (\sum\limits_{i=1}^{a+b}\mu_{i}(T(4;a,b))-(a+b) \overline{d}(T(4;a,b))\right )\\
= & \sqrt{5}(a+b-1)+\sqrt{(a+b)^2-2(a+b)+5}+(a+b)\frac{4}{n}\\&
=\frac{\sqrt{5}}{2}\Big(n-3\Big)+2+\frac{1}{2}\sqrt{n^2-6n+25}-\frac{2}{n}\\&
\ge \frac{\sqrt{5}}{2}\Big(n-3\Big)+2+\frac{1}{2}(n-5)-\frac{2}{n}
>2+\frac{4n}{\pi},
\end{align*}
provided that
\begin{align}\label{x}
\Big(\frac{\sqrt{5}+1}{2}-\frac{4}{\pi}\Big)n> \Big(\frac{3\sqrt{5}+5}{2}+\frac{2}{n}\Big).
\end{align} Since $n\ge 19$ gives that $\frac{2}{n}<0.11$, it follows that the inequality (\ref{x}) holds for all $n\ge 18$. This shows that inequality (\ref{a}) holds for $ T(4;2a,2b)$. If $p=0$, $s_1\ge 2$ and $s_i=1$, for $i=2,3,\dots,r$, then $T$ is the tree $T^{'}$ shown in figure $2$. Applying algorithm $ {\bf (I)} $ to  $ T^{'} $, we find its characteristic polynomial and is given by $
 \phi(T^{'},x)= x(x-1)^{s_1-1}(x^2-3x+1)^{r-2} p(x)$, where $p(x)=(x^4-(r+s_1+5)x^3+(s_1r+4r+3s_1+8)x^2-(2s_1r+5r+2s_1+4)x+s_1+2r).$ Let $\mu_{1}(T^{'})\ge \mu_{2}(T^{'})\cdots\ge \mu_{n}(T^{'})$ be the Laplacian eigenvalues of $T^{'}$.   By Lemma \ref{eigenvalue lower bound}, we have $ \mu_{1}(T^{'})\geq \max\{r,s_1+1\}+1>\overline{d}(T^{'})$ and $ \mu_{2}(T^{'})\geq \min\{r,s_1+1\}>\overline{d}(T^{'})$. Let $x_3\ge x_4$ be the smallest and second smallest zeros of $p(x)$. We have $p(0)=2r+s_1>0$, $p(1)=-s_1(r-1)<0$ and $p(2)=s_1>0$. Therefore, by the intermediate value theorem $x_4\in (0,1)$ and $x_3\in (1,2)$. So, either $\sigma=r$ or $r+1$ for tree $T^{'}$. If $\sigma=r$, then  from the definition of Laplacian energy, we have
 \begin{align*}
 LE(T^{'})=& 2\left (\sum\limits_{i=1}^{\sigma}\mu_{i}(T^{'})-\sigma \overline{d}(T^{'})\right )
 = 2\left (\sum\limits_{i=1}^{r}\mu_{i}(T^{'})-r \overline{d}(T^{'})\right )\\&
 \ge  (\sqrt{5}+1)r+2s_1-(2+2\sqrt{5})+\frac{4r}{n}
 >2+\frac{4n}{\pi},
 \end{align*}
 provided that
 \begin{align}\label{y}
 \Big(\sqrt{5}+1-\frac{8}{\pi}\Big)r+\Big(2-\frac{4}{\pi}\Big)s_1-(4+2\sqrt{5})+\frac{4r}{n}> 0.
 \end{align} It is easy to see that the inequality (\ref{y}) holds for all $r\ge 5$ and $s_1\ge 7$. Since $n=2r+s_1$, it follows that inequality (\ref{y}) holds for all $n\ge 17$.  This shows that the inequality (\ref{a}) holds for $ T^{'}$ in this case. If $\sigma=r+1$, the it can seen by proceeding similarly as above that the  inequality (\ref{a}) holds for $ T^{'}$ in this case as well. If $p=0$, $s_1,s_2\ge 2$ and $s_i=1$, for $i=3,\dots,r$, then $T$ is the tree $T^{''}$ shown in figure $2$. Applying algorithm $ {\bf (I)} $ to  $ T^{''} $, we find its characteristic polynomial and is given by $
  \phi(T^{''},x)= x(x-1)^{s_1+s_2-2}(x^2-3x+1)^{r-3} g(x)$, where $g(x)=(x^6-(r+s_1+s_2+7)x^5+\alpha_1 x^4-\alpha_2 x^3+\alpha_3 x^2+\alpha_4 x+  s_1+s_2+2r-1)$, $\alpha_1=rs_1+rs_2+s_1s_2+5s_1+6r+5s_2+19, \alpha_2=rs_1s_2+4rs_1+3s_1s_2+4rs_2+9s_1+9s_2+14r+24, \alpha_3=2rs_1s_2+5rs_1+3s_1s_2+5rs_2+7s_1+7s_2+16r+13, \alpha_4=2rs_1+2s_1s_2+2rs_2+3s_1+3s_2+9r+1$. Let $\mu_{1}(T^{''})\ge \mu_{2}(T^{''})\dots\ge \mu_{n}(T^{''})$ be the Laplacian eigenvalues of $T^{''}$.  By Lemma \ref{eigenvalue lower bound}, we have
\begin{align*}
  \mu_{1}(T^{''})\geq& \max\{r,s_1+1,s_2+1\}+1>\overline{d}(T^{''}), \mu_{2}(T^{''})\geq \max\{\{r,s_1+1,s_2+1\}\setminus\theta\}>\overline{d}(T^{''})\\&
~~ \text{and}~~   \mu_3(T^{''})\ge \min\{r,s_1+1,s_2+1\}-1>\overline{d}(T^{''}),
\end{align*} where $\theta=\max\{r,s_1+1,s_2+1\}$. Let $x_4\ge x_5\ge x_6$ be the smallest, second smallest and third smallest zeros of $g(x)$. It is easy to see that $x_5,x_6< \overline{d}(T^{''})$, giving that $\sigma=r$ or $r+1$ for tree $T^{''}$. If $\sigma=r$, then  from the definition of Laplacian energy, we have
\begin{align*}
   LE(T^{''})=& 2\left (\sum\limits_{i=1}^{\sigma}\mu_{i}(T^{''})-\sigma \overline{d}(T^{''})\right )
   = 2\left (\sum\limits_{i=1}^{r}\mu_{i}(T^{''})-r \overline{d}(T^{''})\right )\\&
   \ge  (\sqrt{5}+1)r+2s_1+2s_2-( 5+3\sqrt{5})+\frac{4r}{n}
   >2+\frac{4n}{\pi},
   \end{align*}
   provided that
\begin{align}\label{z}
   \Big(\sqrt{5}+1-\frac{8}{\pi}\Big)r+\Big(2-\frac{4}{\pi}\Big)s_1+\Big(2-\frac{4}{\pi}\Big)s_2-(7+3\sqrt{5}-\frac{4}{\pi})+\frac{4r}{n}> 0.
\end{align}
It is easy to see that the inequality (\ref{z}) holds for all $r\ge 2$ and $s_1\ge 8$. Since $n=2r+s_1+s_2-1$, it follows that inequality (\ref{z}) holds for all $n\ge 19$.  This shows that the inequality (\ref{a}) holds for $ T^{''}$. If $\sigma=r+1$, the it can seen by proceeding similarly as above that the  inequality (\ref{a}) holds for $ T^{''}$. This proves the result in this case.   \\
\indent Now, suppose that $T$ is an SNS-tree other than $T(4;2a,2b), T^{'}$ and $T^{''}$. In this case we prove the result by induction on $r$. If $r=2$, then $T$ is a tree shown in figure $3$. Clearly, $T$ has $s=3$ internal vertices, so  using part {\bf (iii)} of Corollary \ref{cor3}, the result holds in this case. Assume that result holds for all trees with $r=k$ vertices of level $1$. We show the result holds for a trees  with $r=k+1$ vertices of level $1$. Let $T$ be a tree with $r=k+1$ vertices (say) $v_1,v_2,\dots,v_k, v_{k+1}$ of level $1$. Let $e=v_0v_{k+1}$ be the edge joining the root $v_0$ with the vertex $v_{k+1}$. Delete the edge $e$, let $T_1$ and $T_2$ be the components of $T-e$. Assume that the order of $T_i$ is $n_i$ and the average vertex degree of $T_i$ is  $\overline{d}(T_i)=2-\frac{2}{n_i}$,  for $i=1,2$, with $n_1\ge n_2$. Since, $ r\ge 3$, therefore $T_1$ is an SNS-tree of order $n_1$ having $r=k$ vertices namely $v_1,v_2,\dots,v_k$ of level $1$, while as $T_2$ is a star $K_{1,s_{k+1}}$. By induction hypothesis the inequality (\ref{a}) holds for both $T_1$ and $T_2$. Let $\sigma_i$ be the number of Laplacian eigenvalues of $T_i$ which are greater than or equal to average degree $\overline{d}(T_i)=2-\frac{2}{n_i}$. Since $T_2$ is a star with at least two vertices and for a star all non-zero Laplacian eigenvalues except the  spectral radius are equal to $1$, it follows that $\sigma_2=1$. To compute $\sigma_1$, we use Applying algorithm $ {\bf (II)} $ with $\alpha=-2+\frac{2}{n_1}$ to the tree $ T_{1}$. Let $u$ corresponds to the pendent vertices in $T_1$. We have
\begin{align*}
a(u)=1-2+\frac{2}{n_1}=-\frac{n_1-2}{n_1}<0,
\end{align*}  as $n_1\ge 4$. This gives that the diagonal entries in the resulting diagonal matrix corresponding each pendent vertex is negative. For the vertices $v_i$, $1\le i\le k$ of level $1$, we have
 \begin{align*}
 a(v_i)=s_i+1-2+\frac{2}{n_1}-\frac{s_i}{a(u)}=2s_i-1+\frac{2s_i}{n_1-2}+\frac{2}{n_1}>0,
 \end{align*} giving that the diagonal entries in the resulting diagonal matrix corresponding each vertex $v_i$ of level $1$ is positive. From Lemma \ref{number1}, it is clear that $\sigma_1=k$ or $k+1$, depending up on whether $a(v_0)<0$ or $a(v_0)>0$. For the root vertex $v_0$, we have
\begin{align}\label{b}
  a(v_0)=&k+p-2+\frac{2}{n_1}-p(-\frac{n_1-2}{n_1})-\sum\limits_{i=1}^{k}\frac{1}{a(v_i)}\\&
  =k+p-2+\frac{2}{n_1}+\frac{p(n_1-2)}{n_1}-\sum\limits_{i=1}^{k}\frac{1}{a(v_i)}.
\end{align}
If $p\ge 1$, then it easy to see that $a(v_0)>0$, giving by Lemma \ref{number1} that $\sigma_1=k+1$. So, suppose that $p=0$ in $T_1$. If at least three $s_i$ say $s_1,s_2,s_3$ are greater or equal $2$, then again it can be seen that $a(v_0)>0$, giving by Lemma \ref{number1} that $\sigma_1=k+1$ in this case as well. So, assume that at most two $s_i$ are greater or equal to $2$, then $T_1\in \{T(4;2a,2b), T^{'},T^{''}\}$, which is not the case.  Thus, it follows that for the tree $T_1$, we have $\sigma_1=k+1$, the number of non-pendent vertices. Let $k_1$ be the number of Laplacian eigenvalues of $T_1$ greater or equal to $\overline{d}(T_1\cup T_2)=2-\frac{4}{n}$. Since $\overline{d}(T_1)\ge \overline{d}(T_1\cup T_2)$, it follows that $k_1\ge \sigma_1$. We claim that $k_1=\sigma_1$. Since $\sigma_1=k+1$, it follows that  $k_1\ge k+1$.  Applying algorithm $ {\bf (II)} $ with $\alpha=-2+\frac{4}{n}$ to the tree $ T_{1}$, we get
\begin{align*}
a(u)=1-2+\frac{4}{n}=\frac{-(n-4)}{n}<0,
\end{align*}
for all $u$ pendent vertices $u$ of $T_1$. This gives that $k_1\le k+1$. Thus we must have $k_1=k+1=\sigma_1$, proving the claim in this case. Further for the tree $T_2$, let $k_2$ be the number of Laplacian eigenvalues greater than or equal to $\overline{d}(T_1\cup T_2)=2-\frac{4}{n}$. Since $\sigma_2=1$ and $T_2$ has at
least one edge, it follows that $k_2=\sigma_2=1$. Thus, we have shown that for the components  $T_1$ and $T_2$ of $T-e$ the inequality (\ref{a}) holds and they satisfies the property that $k_1=\sigma_1$, $k_2=\sigma_2$. Now, applying  Corollary \ref{coru}, it follows that inequality (\ref{a}) holds for $T$, also. Thus, we conclude with help of induction that the result is true for all $r\ge 2$. This completes the proof.\qed

\begin{figure}
\centering
	\includegraphics[scale=.2]{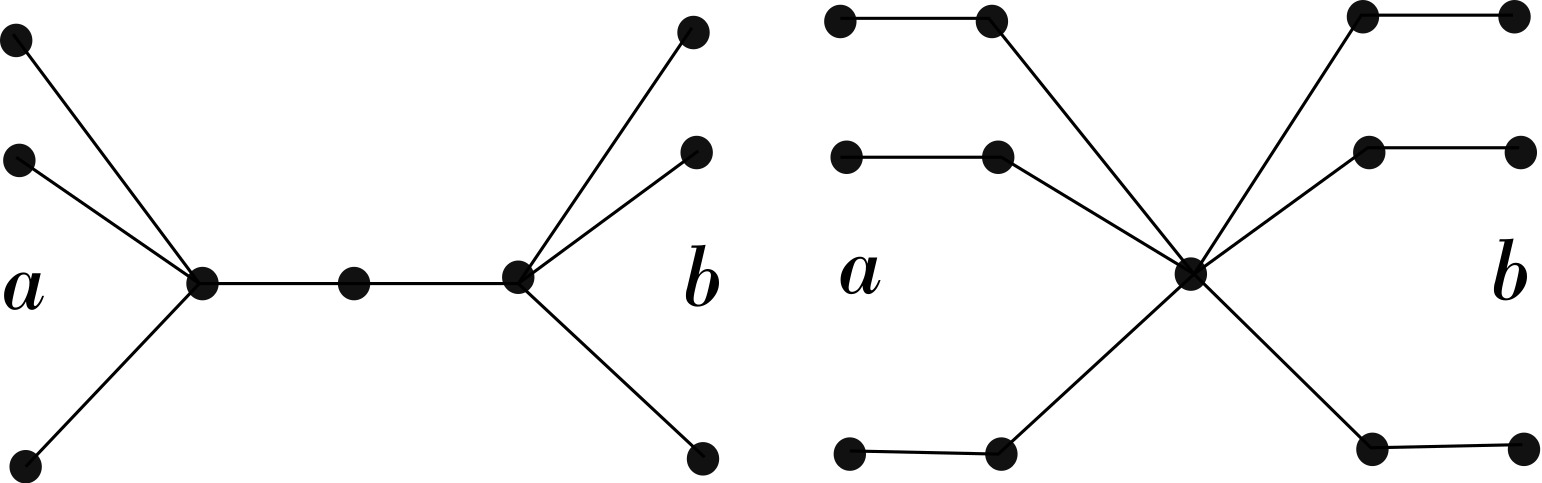}
	\caption{Double broom of diameter $4$ and the tree $T(4;2a,2b)$}
	\label{Figure1}
\end{figure}

\begin{figure}
\centering
	\includegraphics[scale=.2]{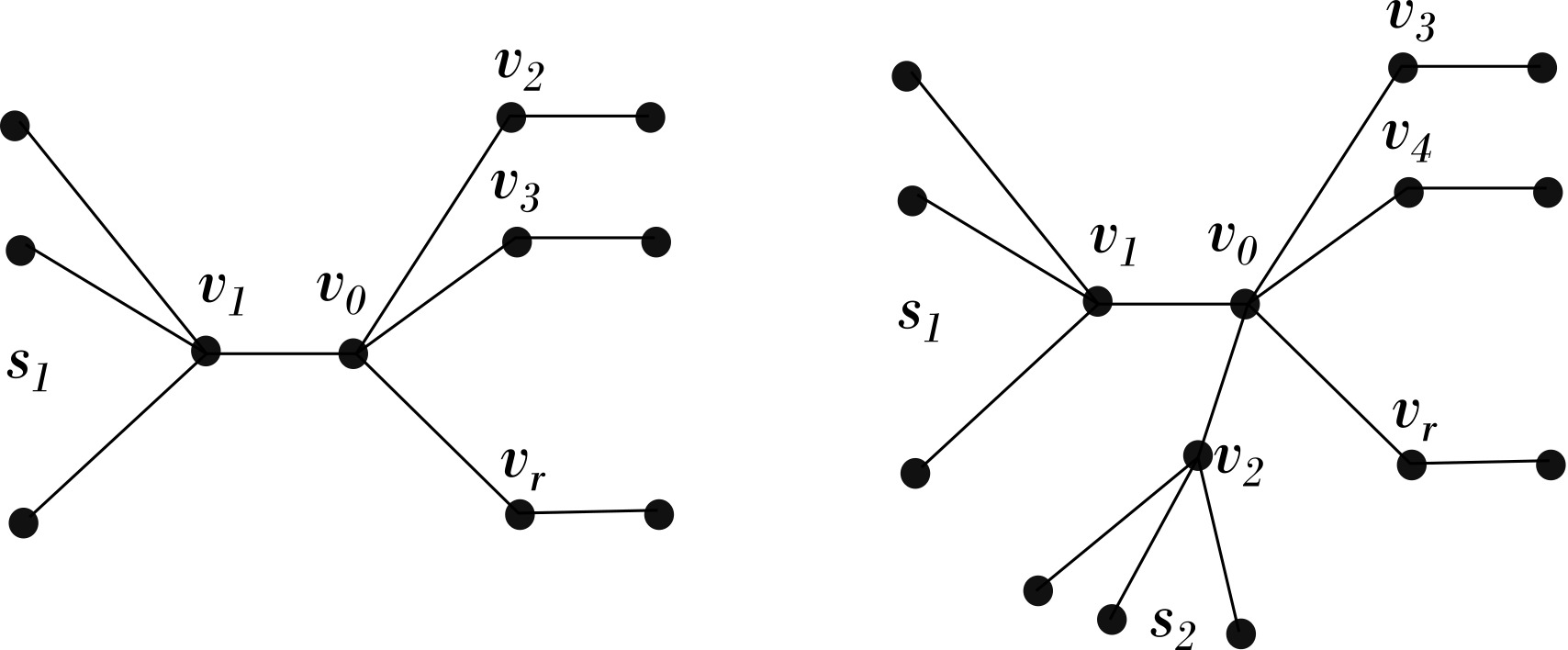}
	\caption{Trees $ T^{'} $ and $ T^{''} $}
	\label{Figure2}
\end{figure}

\begin{figure}
\centering
	\includegraphics[scale=.3]{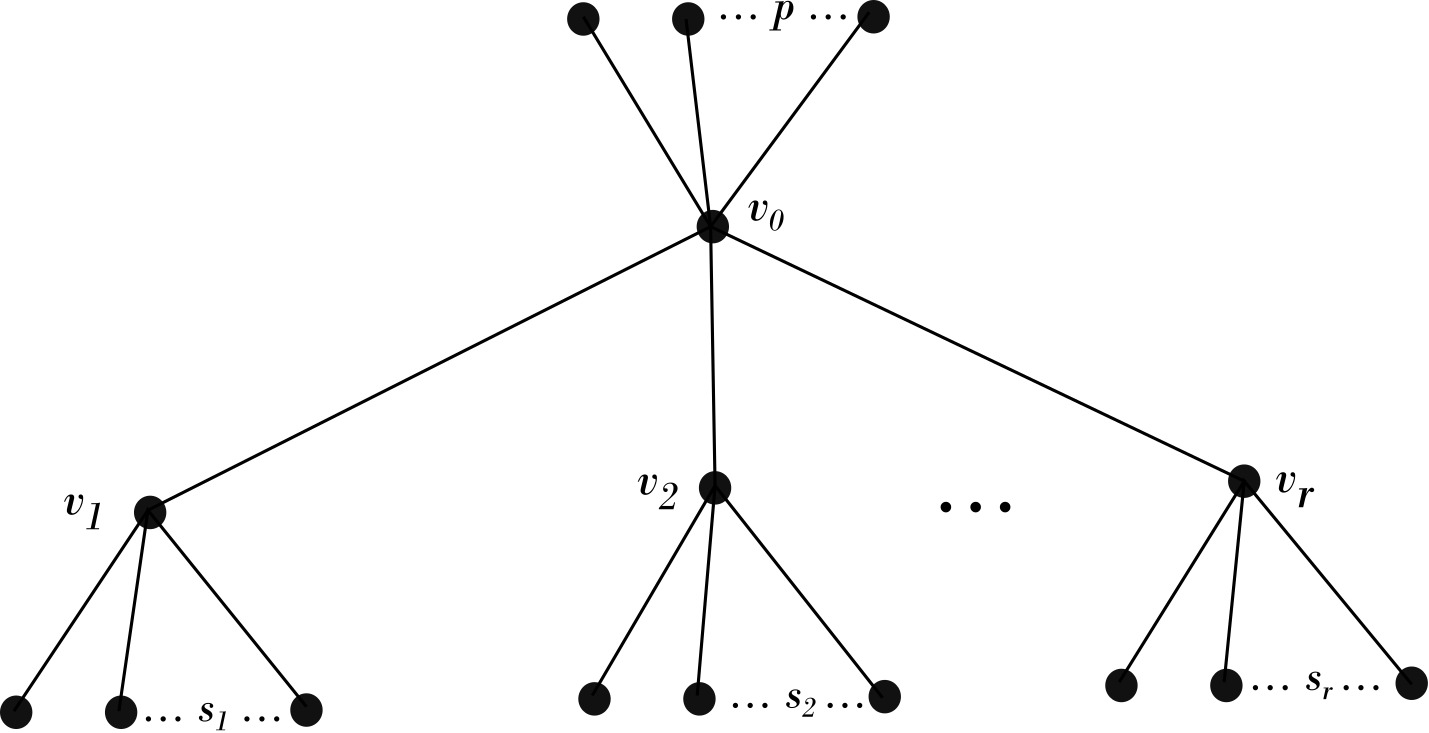}
	\caption{An SNS-tree of diameter $ 4 $}
	\label{Figure3}
\end{figure}

\section{Trees of any diameter}

In this section, we give some sufficient conditions for Conjecture \ref{conjecture} to hold for a tree of order $n$.

The following result gives that if Conjecture \ref{conjecture} holds for two components obtained by deleting a non-pendent edge of a tree $T$ and the number of Laplacian eigenvalues greater than or equal to average degree of the both the components are equal to their number of non-pendent vertices, then Conjecture \ref{conjecture} holds for $T$.
\begin{theorem}\label{thm51}
Let $T_1$ be  a tree of order $n_1$  having $r_1$ non-pendent vertices and let $T_2$ be a tree of diameter at most $3$ having order $n_2$ with $n_1\ge n_2\ge 6$. Let $\sigma_1$ be the number of Laplacian eigenvalues of $T_1$ greater than or equal to average vertex degree $\overline{d}(T_1)=2-\frac{2}{n_1}$. Let $T$ be the tree obtained by joining an edge between a vertex in $T_1$ and a vertex in $T_2$. If $\sigma_1=r_1$, then inequality (\ref{a}) (and hence Conjecture \ref{conjecture}) holds for $T$, provided that $LE(T_1)\ge 2+\frac{4n_1}{\pi}$.
\end{theorem}
\noindent{\bf Proof.} Let $\sigma$ be the number of Laplacian eigenvalues of $T_1\cup T_2$ which are greater than or equal to average vertex degree $\overline{d}(T_1\cup T_2)=2-\frac{4}{n}$, $n_1+n_2=n$. Let $k_i$ be the number of Laplacian eigenvalues of $T_i$ which are greater than or equal to average vertex degree $\overline{d}(T_1\cup T_2)=2-\frac{4}{n}$. Since $n_1\ge n_2$ gives that $\overline{d}(T_1)\ge \overline{d}(T_1\cup T_2)$ and $\overline{d}(T_2)\le \overline{d}(T_1\cup T_2)$, it follows that $k_1\ge \sigma_1$. Let $r_1$ be the number of non-pendent vertices in $T_1$. Suppose that $\sigma_1=r_1$, then $k_1\ge r_1$. We claim that $k_1=\sigma_1=r_1$.  Applying algorithm $ {\bf (II)} $ with $\alpha=-2+\frac{4}{n}$ to the tree $ T_{1}$, we get
\begin{align*}
a(u)=1-2+\frac{4}{n}=\frac{-(n-4}{n}<0,
\end{align*} for all pendent vertices $u$ of $T_1$. This gives that $k_1\le r_1$. Combining $k_1\ge r_1$ and  $k_1\le r_1$, we get $k_1=r_1=\sigma_1$. That proves our claim.  Let $\sigma_2$ be the number of Laplacian eigenvalues of $T_2$ which are greater than or equal to average vertex degree $\overline{d}(T_2)=2-\frac{2}{n_2}$. By given $T_2$ is a tree of diameter at most $3$, it implies that $T_2=K_{1,n_2-1}$ or $T_2=T(a,b)$, a double broom of diameter $3$. If $T_2=K_{1,n_2-1}$, then $\sigma_2=1$, the number of non-pendent vertices. Since  $T_2$ has at
least two edge and $\sigma_2=1$ it follows that $k_2=\sigma_2=1$. On the other hand, if $T_2=T(a,b)$, a double broom of diameter $3$, then as shown in \cite{trevisan}, we have $\sigma_2=2$, the number of non-pendent vertices. Since  $T_2$ is of order $n_2\ge 6$ and $\sigma_2=2$, the number of non-pendent vertices, we must have $k_2=\sigma_2=2$. Thus, we have shown that the trees  $T_1$ and $T_2$ have the  property that $k_1=\sigma_1$ and  $k_2=\sigma_2$. Since the inequality (\ref{a}) always holds for $T_2$, by Theorem \ref{thm23}, it follows by using  Corollary \ref{coru} that inequality (\ref{a}) (and hence Conjecture \ref{conjecture})  holds for $T$, provided that $LE(T_1)\ge 2+\frac{4n_1}{\pi}$. This completes the proof.\qed

If $T_2$ is an SNS-tree of diameter $4$ other than the trees $T(4;2a,2b), T^{'},T^{''}$ defined in section $4$, the using Theorem \ref{thm42}, and proceeding similarly as in Theorem \ref{thm51}, we arrive at the following result.

\begin{theorem}\label{thm52}
Let $T_1$ be  a tree of order $n_1$  having $r_1$ non-pendent vertices and let $T_2$ be an SNS-tree of diameter $4$ other than the trees $T(4;2a,2b), T^{'}, T^{''}$ having order $n_2$ with $n_1\ge n_2\ge 6$. Let $\sigma_1$ be the number of Laplacian eigenvalues of $T_1$ greater than or equal to average vertex degree $\overline{d}(T_1)=2-\frac{2}{n_1}$. Let $T$ be the tree obtained by joining an edge between a vertex in $T_1$ and a vertex in $T_2$. If $\sigma_1=r_1$, then inequality (\ref{a}) (and hence Conjecture \ref{conjecture}) holds for $T$, provided that $LE(T_1)\ge 2+\frac{4n_1}{\pi}$.
\end{theorem}

The following result gives that if Conjecture \ref{conjecture} holds for two components obtained by deleting a non-pendent edge of a tree $T$ and one of the components $T_1$ satisfy $\mu_{\sigma_1+1}(T_1)- \overline{d}(T_1)<\frac{-2}{n}$ and other component $T_2$ have the property that $\sigma_2$ is  equal to the number of non-pendent vertices, then Conjecture \ref{conjecture} holds for $T$. Note that $\sigma_i$ is the  number of Laplacian eigenvalues of $T_i$ greater than or equal to average vertex degree $\overline{d}(T_i)$.
\begin{theorem}
Let $T_1$ be  a tree of order $n_1$  having Laplacian eigenvalues $\mu_1(T_1)\ge\dots\ge \mu_{n_1-1}(T_1),\mu_{n_1}(T_1)=0$ and let $T_2$ be a tree of diameter at most $3$ or an SNS-tree of diameter $4$ other than the trees $T(4;2a,2b), T^{'},T^{''}$  having order $n_2$ with $n_1\ge n_2\ge 6$. Let $\sigma_1$ be the number of Laplacian eigenvalues of $T_1$ greater than or equal to average vertex degree $\overline{d}(T_1)=2-\frac{2}{n_1}$.
 Let $T$ be the tree of order $n=n_1+n_2$ obtained by joining an edge between a vertex in $T_1$ and a vertex in $T_2$. If $\mu_{\sigma_1+1}- \overline{d}(T_1)<\frac{-2}{n}$, then inequality (\ref{a}) (and hence Conjecture \ref{conjecture}) holds for $T$, provided that $LE(T_1)\ge 2+\frac{4n_1}{\pi}$.
\end{theorem}
{\bf Proof.} Let $\sigma$ be the number of Laplacian eigenvalues of $T_1\cup T_2$ which are greater than or equal to average vertex degree $\overline{d}(T_1\cup T_2)=2-\frac{4}{n}$, $n_1+n_2=n$. Let $k_i$ be the number of Laplacian eigenvalues of $T_i$ which are greater than or equal to average vertex degree $\overline{d}(T_1\cup T_2)=2-\frac{4}{n}$. Since $n_1\ge n_2$ gives that $\overline{d}(T_1)\ge \overline{d}(T_1\cup T_2)$ and $\overline{d}(T_2)\le \overline{d}(T_1\cup T_2)$, it follows that $k_1\ge \sigma_1$. Since, $\mu_{\sigma_1+1}- \overline{d}(T_1)<\frac{-2}{n}$ gives that $\mu_{\sigma_1+1}- \overline{d}(T_1\cup T_2)=\mu_{\sigma_1+1}-2+\frac{4}{n}=\mu_{\sigma_1+1}-2+\frac{2}{n_1}+\frac{4}{n}-\frac{2}{n_1}<\frac{-2}{n}+\frac{4}{n}-\frac{2}{n_1}<0$, as $n>n_1$. It follows that $k_1\le \sigma_1$ and so we must have $\sigma_1=k_1$. Now, proceeding similar to Theorem \ref{thm51} the result follows.\\

\noindent{\bf Acknowledgements.} The research of S. Pirzada is supported by SERB-DST, New Delhi under the research project number  MTR/2017/000084.

\end{document}